\documentclass[12pt,notitlepage]{amsart}
\usepackage{latexsym,amsfonts,amssymb,amsmath,amsthm}
\usepackage{color}
\usepackage{hyperref}
\hypersetup{colorlinks=true,citecolor=blue,linkcolor=blue,urlcolor=blue,
pdfstartview=FitH}
\usepackage[inner=1.0in,outer=1.0in,bottom=1.0in, top=1.0in]{geometry}

\newcommand{\mz}{\ensuremath{\mathbb Z}}
\newcommand{\mr}{\ensuremath{\mathbb R}}

\newcommand{\onehalf}{\ensuremath{ \frac{1}{2}}}

\newcommand{\SL}[2]{SL_{#1}(\mathbb{#2})}

\newcommand{\sgn}{\operatorname{sgn}}

\theoremstyle{plain}		
	\newtheorem{mytheo}{Theorem}[section]
	
	\newtheorem{myprop}[mytheo]{Proposition}
	
     \newtheorem{mylemma}[mytheo]{Lemma}

\theoremstyle{remark}

\newcommand{\tn}{\ensuremath{|\theta N|}}
\numberwithin{equation}{section}
\begin{document}
\author{Daniel Godber}
\address{Department of Mathematics, Texas A\&{}M University, College Station, TX 77843-3368, U.S.A.}
\email{dgodber@math.tamu.edu}

\thanks{}

\begin{abstract}
We study sums of the form $\sum_{n\le N} a(n)e^{2\pi i\alpha n}$, where $\alpha$ is any real number and the $a(n)$ are the Fourier coefficients of either a holomorphic cusp form, a Maass cusp form, or the symmetric-square lift of a holomorphic cusp form. We obtain bounds that are uniform in both $\alpha$ and the form itself. We also improve a bound on a sum of the form $\sum_{n\le N} a(n)e^{2\pi i(\alpha n+\beta n^\theta)}$, where the $a(n)$ are the Fourier coefficients of a holomorphic cusp form, $\alpha$ and $\beta$ are any real numbers, and $0 \le \theta < 1$. This last bound is uniform in $\alpha$, but not with respect to the form.
\end{abstract}

\title{Additive Twists of Fourier Coefficients of Modular Forms}
\maketitle

\section{Introduction}
\noindent Let $f$ be a holomorphic Hecke eigenform of weight $k$ for $SL(2,\mz)$ with Fourier expansion
\begin{align}
\label{eq:firstf}
f(z) = \displaystyle\sum_{n=1}^{\infty}\lambda_f(n) n^{\frac{k-1}{2}}e(nz), \text{ } e(z) = e^{2\pi iz}
\end{align}
(for Im$z > 0$), normalized so that $\lambda_f(1)=1$. It is a well-known result that for any real number $\alpha$,
\begin{align}
\label{eq:knownbound}
 \sum_{n\le N} \lambda_f(n) e(\alpha n) \ll_{f} N^{1/2}\log 2N.
\end{align}
 The uniformity in $\alpha$ is convenient because it allows us to obtain the same bound for a sum of Fourier coefficients restricted to any arithmetic progression. The uniformity also suggests that there is no correlation between the Fourier coefficients and additive characters. It would be interesting and useful for certain applications to make the dependence on the form explicit.

While it is simple to prove \eqref{eq:knownbound} by using the ideas in \cite[Theorem 5.3]{Iw1}, the method unfortunately cannot be modified to obtain a bound that is uniform with respect to the form. The proof of \eqref{eq:knownbound} essentially relies only on an estimate of the size of $f(z)$ and partial summation. To obtain uniformity with respect to the form, we require some heavier machinery such as the Voronoi summation formula and careful analysis of certain exponential integrals.

\begin{mytheo}
\label{theorem:first}
Let $f(z)$ be as in \eqref{eq:firstf}. Then for any real number $\alpha$ and any $\varepsilon > 0$,
\begin{align}
\label{eq:gl2holosum}
\displaystyle\sum_{n\le N} \lambda_f(n) e(\alpha n) \ll_\varepsilon N^{1/2 + \varepsilon}\mathfrak{q}_f^{1/4+\varepsilon},
\end{align}
where $\mathfrak{q}_f\asymp k^2$ is the analytic conductor for $f$.
\end{mytheo}
\noindent A similar statement holds for Maass forms.
\begin{mytheo}
\label{theorem:second}
Let $u(z)$ be a Hecke-Maass form with Laplace eigenvalue $\frac{1}{4} +T^2$, and with Hecke eigenvalues $\lambda_u(n)$. Then if $u$ has the expansion
\begin{align}
\label{eq:massexpansion}
u(z) = y^{1/2}\sum_{n\ne 0}\lambda_u (n)K_{iT}(2\pi |n|y)e(nx),
\end{align}
then for any real number $\alpha$ and any $\varepsilon > 0$,
\begin{align}
\displaystyle\sum_{n\le N} \lambda_u(n) e(\alpha n) \ll_\varepsilon N^{1/2 + \varepsilon}\mathfrak{q}_u^{1/4+\varepsilon},
\end{align}
where $\mathfrak{q}_u\asymp T^2$ is the analytic conductor for $u$.
\end{mytheo}
\noindent This is an improvement over the bound $N^{1/2+\varepsilon}\mathfrak{q}_u^{1/2}$ which appears in \cite[\S 8.3]{Iw2}.

It is known due to Gelbart and Jacquet \cite{GJ} that the symmetric-square lift $L(F,s) := L(\text{Sym}^2 f,s)$ is also an L-function for some $GL(3,\mz)$-automorphic representation $F$. This L-function satisfies the functional equation
\begin{align}
\label{eq:functionalequation}
\Lambda_F(s) :=  \pi^{-\frac{3s}{2}}\Gamma(\frac{s+1}{2})\Gamma(\frac{s+k-1}{2})\Gamma(\frac{s+k}{2})L(F,s) = \Lambda_F(1-s),
\end{align} 
and its coefficients are given by $A_F(1,n) = \sum_{ml^2 = n}\lambda_f(m^2)$.

\begin{mytheo}
\label{theorem:third}
Let $F$ be the symmetric-square lift of a holomorphic Hecke eigenform of weight $k$ for $\SL{2}{Z}$ with $A_F(1,1) = 1$. Then for any real number $\alpha$ and any $\varepsilon>0$, 
\begin{align}
\label{eq:mainresult}
\displaystyle\sum_{n\le N} A_F(1,n) e(\alpha n) \ll_\varepsilon N^{3/4 + \varepsilon}\mathfrak{q}_F^{1/4+\varepsilon},
\end{align}
where $\mathfrak{q}_F\asymp k^2$ is the analytic conductor for $F$.
\end{mytheo}

A similar result for when $F$ is the lift of a Maass form was proven by Xiaoqing Li and M. Young \cite{LY1}, and our proof will closely follow theirs. They obtained the bound $N^{3/4+\varepsilon}\mathfrak{q}_F^{D+\varepsilon}$, where $D = 1/4$ assuming the Ramanujan conjecture, and $D = 1/3$ unconditionally. Our result is stronger because in our case we can use Deligne's bound. Prior to \cite{LY1}, Miller \cite{M} had obtained the bound $N^{3/4+\varepsilon}$, where the bound is uniform in $\alpha$ but has an implied constant that depends on $F$. Xiannan Li \cite{Li} generalized the result of \cite{LY1} by considering the case when $F$ is a general $GL(3,\mz)\backslash GL(3,\mr)$ cusp form, obtaining the bound $N^{3/4+\varepsilon}\mathfrak{q}_F^D$, where $D=1/4$ assuming the Ramanujan conjecture, and $D=5/12$ unconditionally.

Li and Young were motivated to investigate the non--holomorphic $GL(3)$ case by a previous paper \cite{LY2}, which required the application of the $GL(3)$ Voronoi formula with a varying form. Until now, the $GL(2)$ case had been previously unstudied and as it turns out, it is nontrivial.

The main tool in the proofs of Theorems \ref{theorem:first}-\ref{theorem:third} is the Voronoi summation formula, which relates the sum of the Hecke eigenvalues of our forms to another sum with a weight function given as an integral transform (see for example Theorem \ref{theorem:gl2voronoiholo}). This integral contains a ratio of gamma factors that can be estimated by Stirling's approximation, as well as a Mellin transform that can be estimated by the method of stationary phase. These estimations leave us with an exponential integral of the form $\int_\alpha^\beta g(t)e^{ih(t)}dt$.

To estimate such an integral, we use two well-known lemmas.  The first is Lemma 5.1.2 \cite{H}, which states that if $|h'(x)| \gg \kappa > 0$ on $[\alpha, \beta]$, then $\int_\alpha^\beta g(t)e^{ih(t)}dt \ll V/\kappa$ where $V$ is the total variation of $g$ on $[\alpha, \beta]$ plus the maximum modulus of $g$ on $[\alpha, \beta]$. The second is Lemma 5.1.3 \cite{H}, which states that if $|h''(x)| \gg \lambda > 0$ on $[\alpha, \beta]$, then $\int_\alpha^\beta g(t)e^{ih(t)}dt \ll V\lambda^{-\onehalf}$, where $V$ is as before.

In the case of Theorem \ref{theorem:second}, we can directly apply these lemmas to bound the exponential integral. In the proof of Theorem \ref{theorem:third}, the lemmas are not quite good enough by themselves. However, we can exploit the fact that the undesirable bounds occur for only a short interval of the summation to obtain our result. It is convenient that there is no short interval behavior in the case of Theorem \ref{theorem:second}, so that we do not have to assume the Ramanujan conjecture to obtain our result.

The exponential integral that arises in the proof of Theorem \ref{theorem:first} requires much more work because the phase function $h(t)$ may have a degenerative stationary point. That is, the first and second derivatives of the phase function may have a common zero (or two nearby zeros), making the above lemmas alone incapable of obtaining \eqref{eq:gl2holosum}. In order to estimate this integral, we use a modified method of stationary phase, out of which the Airy function naturally appears (see  \eqref{eq:airy} below). Using the properties of the Airy function and again exploiting the fact that some of the bounds obtained occur for only a short interval of the summation, we arrive at our result.\footnote{Nicolas Templier has informed the author of an alternate proof of Theorem \ref{theorem:first} that uses the bound $\|y^{k/2}f\|_\infty \ll k^{1/4+\varepsilon}\|f\|_2$ due to Xia \cite{Xi}. However, the methods used in our proof are important for understanding the behavior of the Voronoi summation formula when the underlying form is varying. Moreover, while the final results stated in our theorem are uniform in $\alpha$, in the course of our proof we obtain bounds for specific rational approximations of $\alpha$.}

It is surprising that the result for $GL(2)$ holomorphic forms (Theorem \ref{theorem:first}) is far more difficult to obtain than the corresponding result for non--holomorphic Maass forms (Theorem \ref{theorem:second}). In fact the result even requires more work than the $GL(3)$ case.

We end the paper with a short note on how to improve a bound given by Sun \cite{Su} where the sum in Theorem \ref{theorem:first} is twisted with a nonlinear exponential term.
\begin{mytheo}
\label{theorem:sun}
Fix $0\le\theta<1$. Then for any $\alpha,\beta \in \mr$ and $N$ sufficiently large, we have that
\begin{align}
\label{eq:sunbound}
\displaystyle\sum_{n\le N} \lambda_f(n) e(\beta n^\theta + \alpha n) \ll N^{1/2 + \theta/2 +\varepsilon}, 
\end{align}
where the implied constant depends only on $\beta,\theta, \varepsilon$ and $f$.
\end{mytheo}
This is an improvement over Sun's previous bound $N^{1-\theta/2+\varepsilon}$, valid for only $0 < \theta \le \onehalf$, which was obtained by a more complicated method. Liu and Ren \cite{LR} have noted in passing that Sun's bound could be improved to $N^{1/2+\theta+\varepsilon}$ when $\theta \le 1/3$ by a simple application of partial summation. Note that this bound is not uniform in the weight of the form. In fact the proof is completely different than those of the previous theorems in this paper. The proof relies on a convenient bound for the size of $f(z)$ (see \eqref{eq:fbound} below) and the estimation of certain exponential sums, much like the proof of \eqref{eq:knownbound} found in \cite[\S 5.1]{Iw1}. We will actually prove a more general statement (see Theorem \ref{theorem:sungeneral} below) from which Theorem \ref{theorem:sun} immediately follows.

By Iwaniec-Luo-Sarnak \cite[Appendix C]{ILS}, the bound \eqref{eq:sunbound} is essentially sharp for $\theta = \onehalf$. More precisely, if  $w(x)$ is a fixed smooth weight function compactly supported on $\mr^+$, then
\begin{align}
\displaystyle\sum_{n\ge 1} \lambda_f(n) e(-2\sqrt{n})w(\frac{n}{N}) =  CN^{3/4} + O(N^{1/4+\varepsilon})
\end{align}
for some constant $C$.

Theorems \ref{theorem:sun} and \ref{theorem:sungeneral} can actually be generalized in the following way. Let $\Gamma$ be a general discrete group for which $\infty$ is a cusp of width 1, and $\vartheta$ be a multiplier system of weight $k>0$ for $\Gamma$ that is singular at the cusp $\infty$. (See \cite[\S 2.3, 2.6 and 5.1]{Iw1} for definitions and details). Then if $f(z)$ is a cusp form for $\Gamma$ with respect to the multiplier system $\vartheta$, then \eqref{eq:fbound} and hence Theorems \ref{theorem:sun} and \ref{theorem:sungeneral} still hold true.

\section{Acknowledgements}
\noindent The author would like to thank Matthew P. Young for all of his valuable suggestions and encouragement.

\section{Proof of Theorem 1.1}
\noindent In this section we will prove the following result from which Theorem \ref{theorem:first} can be easily deduced by an unsmoothing argument (see Lemma 9 \cite{LY1}):
\begin{mytheo}
\label{theorem:firstsmoothed}
Let $f$ be as in \eqref{eq:firstf} and let $w$ be a weight function satisfying 
\begin{align}
\begin{cases}
\label{eq:wconditions}
\text{w is smooth with compact support on } [N,2N],\\
|w^{(j)}(y)| \le c_jN^{-j},
\end{cases}
\end{align}
for all $j = 0,1,2,\ldots$, where the $c_j$ are some positive real numbers. Then for any real number $\alpha$ and any $\varepsilon>0$,
\begin{align}
\label{eq:smoothedresult}
\displaystyle\sum_{n\ge 1} \lambda_f(n)e(\alpha n)w(n) \ll_{\varepsilon,c_j} N^{1/2 + \varepsilon}\mathfrak{q}_f^{1/4+\varepsilon},
\end{align}
where $\mathfrak{q}_f\asymp k^2$ is the analytic conductor for $f$.
\end{mytheo}

\subsection{$GL(2)$ Voronoi formula}
We will need the following version of the $GL(2)$ Voronoi formula, which is equivalent to the more familiar formula involving the Bessel function.
\begin{mytheo} \cite[Equations (1.12), (1.15)]{MS2}
\label{theorem:gl2voronoiholo}
Let $\psi(x)$ be a smooth function with compact support on the positive reals. Let $d,\overline{d},c \in \mz$ with $c\ne0, (c,d) = 1,$ and $d\overline{d}\equiv1 \pmod{c}$. Then 
\begin{align}
\displaystyle\sum_{n\ge 1}\lambda_f(n)e(\frac{n\overline{d}}{c}) \psi(n) = c \displaystyle\sum_{n\ge 1}\frac{\lambda_f(n)}{n}e(-\frac{nd}{c}) \Psi(\frac{n}{c^2}),
\end{align}
where for $\sigma > -1-(k+1)/2$
\begin{align}
\label{eq:psifirst}
\Psi(x) =  i^{k-1}\frac{1}{2\pi^2}\displaystyle\int_{(\sigma)}(\pi^2x)^{-s}\frac{\Gamma(\frac{1+s+(k+1)/2}{2})\Gamma(\frac{1+s+(k-1)/2}{2})}{\Gamma(\frac{-s+(k+1)/2}{2})\Gamma(\frac{-s+(k-1)/2}{2})}\tilde{\psi}(-s)ds,
\end{align}
where $\tilde{\psi}(s)$ is the Mellin transform of $\psi(x)$. 
\end{mytheo}

In truth, directly applying Miller and Schmid's formula will give what at first appears to be a different formula for $\Psi$. In particular, the gamma factors in \eqref{eq:psifirst} are different. However, using the relation
\begin{align}
\label{eq:gammarelation}
\Gamma(1-s)\Gamma(s) = \pi\csc\pi s
\end{align}
and the fact that $k$ is even, one can rewrite the gamma factors and restate the formula as given above. Writing the formula as above has two advantages. First, since the arguments of the gamma factors all lie in the right half-plane, we can easily apply Stirling's formula. Second, the relation between the gamma factors in the integral and the gamma factors in the functional equation for the L-function of a holomorphic form is much more obvious. Recall that the functional equation for $L(f,s)=\displaystyle\sum_{n\ge 1}\frac{\lambda_f(n)}{n^s}$ is given by 
\begin{align}
\Lambda_f(s) := \pi^{-s}\Gamma(\frac{s+(k-1)/2}{2})\Gamma(\frac{s+(k+1)/2}{2})  =  \Lambda_f(1-s).
\end{align}



%
%
%
\subsection{Bounds on the Gamma Factors}
For the benefit of the reader, we explicitly calculate the asymptotic expansions of the gamma factors using Stirling's formula. Write $s=\sigma -i\tau$.
Fix $\sigma > -1$ and let $C$ be a nonnegative real number large enough so that $1+\sigma+C-\delta>0$ and $-\sigma+C-\delta>0$, where $\delta$ is some fixed real number. The purpose of this number is to allow for some small variations of $C$. For example, to calculate the gamma factors containing $k$ and $k-1$, we can use $C = k$ in both cases and adjust $\delta$ by 1. Now when at least one of $|\tau|$ or $C$ is approaching $\infty$, we can apply Stirling's formula to get

\begin{align}
\log\Gamma(\frac{1+\sigma-i\tau+C-\delta}{2}) & = \frac{\sigma - i\tau + C -\delta}{2}\log(\frac{1+\sigma-i\tau+C-\delta}{2})\\\nonumber &+\frac{-1-\sigma+i\tau-C+\delta}{2} + \onehalf\log2\pi + \sum_{j=1}^{M-1} \frac{c_j}{(C-i\tau)^{j}} +  O(\frac{1}{|C-i\tau|^{M}})
\end{align}
for some constants $c_j$. Now notice that 
\begin{align}
\log(\frac{1+\sigma-i\tau+C-\delta}{2}) & = \log(\frac{C-i\tau}{2}) + \log(1+\frac{1+\sigma-\delta}{C-i\tau}) \\\nonumber
&= \log(\frac{C-i\tau}{2}) + \frac{1+\sigma-\delta}{C-i\tau} +  \sum_{j=2}^{M-1} \frac{d_j}{(C-i\tau)^{j}} +  O(\frac{1}{|C-i\tau|^{M}})
\end{align}
Hence we can write
\begin{align}
\label{eq:loggamma1}
\log\Gamma(\frac{1+\sigma-i\tau+C-\delta}{2}) & = \frac{\sigma - i\tau + C -\delta}{2}\log(\frac{C-i\tau}{2}) +\frac{-C+i\tau}{2}\\\nonumber &+ \onehalf\log2\pi + \sum_{j=1}^{M-1} \frac{C_j}{(C-i\tau)^{j}} +  O(\frac{1}{|C-i\tau|^{M}})
\end{align}
for some constants $C_j$. By exponentiating \eqref{eq:loggamma1}, we can calculate that
\begin{align}
\label{eq:stirlingmain}
\frac{\Gamma(\frac{1+\sigma-i\tau+C-\delta}{2})}{\Gamma(\frac{-\sigma+i\tau+C-\delta}{2})}  = \left(\frac{\sqrt{C^2+\tau^2}}{2}\right)^{\sigma + \onehalf}&e^{-i\tau\log\frac{\sqrt{C^2+\tau^2}}{2e}}e^{i(\delta+\onehalf-C)\arctan(\frac{\tau}{C})}\\\nonumber
&\times\left(c + \frac{P_1(C,\tau)}{C^2+\tau^2} + \frac{P_2(C,\tau)}{(C^2+\tau^2)^2} + \cdots + O(\max(C,|\tau|)^{-A})\right),
\end{align}
where $c$ is some constant and each $P_j(C,\tau)$ is a polynomial of degree $j$. The constant $c$ and the polynomials depend only on $\delta$. 

We note down here an asymptotic expansion for \eqref{eq:psifirst}:
\begin{align}
\label{eq:psifirstasymp}
\Psi(x) = \frac{i^{k-1}}{2\pi^2} \displaystyle\int_{(\sigma)}(\pi^2x)^{-s} &\left( \frac{\tau^2+(k/2)^2}{2^2}\right)^{\sigma+\onehalf} e^{-i\tau\log(\frac{\tau^2+(k/2)^2}{(2e)^2})} e^{i(1-k)\arctan(\frac{2\tau}{k})} \tilde{\psi}(-s)ds\\\nonumber
&\times\left(c + Q_1(k,\tau) + Q_2(k,\tau) + \cdots + Q_{A-1}(k,\tau)  + O(\max(k,|\tau|)^{-A})\right).
\end{align}
where $c$ is some absolute constant and each $Q_j(k,\tau) = O(\max(k,|\tau|)^{-j})$ is a rational function. 

\subsection{Bounding $\Psi(x)$ and $S$}

Let $Q \ge 1$ be a parameter to be chosen later. By Dirichlet's approximation theorem, there exist coprime integers $a,q$ with $1\le q \le Q$ such that $\alpha = \frac{a}{q} + \frac{\theta}{2\pi}$ with $|\frac{\theta}{2\pi}| \le (qQ)^{-1}$. Then we can rewrite the left-hand side of \eqref{eq:smoothedresult}
\begin{align}
S = \displaystyle\sum_{n\ge 1} \lambda_f(n) e(\frac{an}{q}) \psi(n),
\end{align}
where 
\begin{align}
\label{eq:psi2}
\psi(y) = e^{i\theta y}w(y).
\end{align} 


\begin{mylemma}
\label{lemma:gl2bound}
Let $\psi(x)$ be defined by \eqref{eq:psi2} and define
\begin{align}
U = \max(k^2, \tn^2)
\end{align}
and
\begin{align}
\Delta& = \left |xN-\frac{1}{(2\pi)^2}\tn k \right |.
\end{align}
Then 
\begin{align}
\Psi(x) \ll \mathcal{M + E},
\end{align}
where
\begin{align}
\label{eq:gl2boundmain}
\mathcal{M} = U^{1/2}|Nk|^\varepsilon(1+\frac{xN}{U(Nk)^\varepsilon})^{-A},
\end{align}
and $\mathcal{E}=0$ unless $k^{1-\varepsilon} \le \tn \le k^{1+\varepsilon}$, in which case
\begin{align}
\mathcal{E} = 
\begin{cases}
\label{eq:firstcases}
k^{7/6+\varepsilon} & \text{if }  \Delta \ll k^{4/3+\varepsilon},\\
\frac{k^{3/2+\varepsilon}}{\Delta^{1/4}} & \text{if } k^{4/3+\varepsilon} \ll \Delta \ll k^{2+\varepsilon},\\
0 & \text{otherwise}.
\end{cases}
\end{align}
\end{mylemma}
\noindent Deferring the proof of Lemma \ref{lemma:gl2bound}, we first prove Theorem \ref{theorem:firstsmoothed}.\\
{\it Proof of Theorem \ref{theorem:firstsmoothed}.}
  By Theorem \ref{theorem:gl2voronoiholo} and Lemma \ref{lemma:gl2bound}, we have that $S \ll S_\mathcal{M} + S_{\mathcal{E}}$, corresponding to $\Psi \ll \mathcal{M} + \mathcal{E}$. It is easy to see that 
\begin{align}
S_\mathcal{M} \ll q(k + \tn)(NkQ)^\varepsilon,
\end{align}
and since $q|\theta| \le 2\pi Q^{-1}$, we have that
\begin{align}
S_\mathcal{M} \ll (Qk + Q^{-1}N)(NkQ)^\varepsilon.
\end{align}
It remains to bound $S_{\mathcal{E}}$. Now in the case that $\Delta \ll k^{4/3+\varepsilon}$, applying Deligne's bound gives us the bound
\begin{align}
qk^{7/6}k^{-2/3}(Nkq)^\varepsilon \ll qk^{1/2}(Nkq)^\varepsilon.
\end{align}
In the case that $k^{4/3+\varepsilon} \ll \Delta \ll k^{2+\varepsilon}$, we assume that $\Delta \in [Y,2Y]$ and divide this interval into $\ll Yk^{-4/3-\varepsilon}$ subintervals of length at most $k^{4/3+\varepsilon}$. Then applying Deligne's bound again gives us $O(Y k^{-4/3-\varepsilon})$ instances of bounds of the form
\begin{align}
q\frac{k^{3/2}}{Y^{1/4}}k^{-2/3}(Nkq)^\varepsilon
\end{align}
so that the sum of the bounds is bounded by $qk(Nkq)^\varepsilon$. Putting this together we have that 
\begin{align}
S \ll (Qk + Q^{-1}N)(Nkq)^\varepsilon.
\end{align}
Choosing $Q = N^{1/2}k^{-1/2}$ gives the bound stated in Theorem \ref{theorem:firstsmoothed}. \qed\\ 

For the proof of Lemma \ref{lemma:gl2bound}, we will require the following additional lemma, which can be proven by the method of stationary phase. For a proof see \cite[Lemma 5.1]{LY1}.

\begin{mylemma}
\label{lemma:mellin}
Let $\tau, \theta$ and $N$ be real numbers and let $w$ be as in Theorem \ref{theorem:firstsmoothed}. Let
\begin{align}
I = \displaystyle\int_0^\infty w(x)e^{i\theta x}x^{i\tau} \frac{dx}{x}.
\end{align}
If $|\tau| \ge 1$ and $\tn \ge 1$, then 
\begin{align}
\label{eq:mellinasymp}
I = \sqrt{2\pi} w(-\tau/\theta)|\tau|^{-\onehalf} e^{i\tau\log |\tau/(e\theta)|}e^{\frac{i\pi}{4}\sgn(\theta)} + O(|\tau|^{-3/2}).
\end{align}
Furthermore, if $|\tau| \ge \tn^{1+\varepsilon}$ then 
\begin{align}
I \ll_{A,\varepsilon} |\tau|^{-A}
\end{align}
and if $|\tau| \le \tn^{1-\varepsilon}$ then
\begin{align}
I \ll_{A,\varepsilon} \tn^{-A}.
\end{align}
\end{mylemma}

Note that if $\tn \le 1$ then $w_\theta(x) := w(x)e^{i\theta x}$ satisfies the same properties as $w(x)$ and so $I = \tilde{w_\theta}(i\tau)$. Integrating by parts shows that
\begin{align}
I \ll_A (1 + |\tau|)^{-A}.
\end{align}

\noindent{\it Proof of Lemma \ref{lemma:gl2bound}.} First note that $\tilde{\psi}(-\sigma +i\tau) = \displaystyle\int_0^\infty w(x)x^{-\sigma}e^{i\theta x}x^{i\tau}\frac{dx}{x}$.
If $\tn \le 1$ then by the modifying the remark after Lemma \ref{lemma:mellin}, we see that
\begin{align}
\tilde{\psi}(-\sigma +i\tau) \ll_{A,\sigma} N^{-\sigma}(1+|\tau|)^{-A}.
\end{align}
If $\tn > 1$, then we can apply Lemma \ref{lemma:mellin}. To unify all the cases we use the bound 
\begin{align}
\tilde{\psi}(-\sigma +i\tau) \ll_{A,\varepsilon, \sigma} N^{-\sigma}(1+\frac{|\tau|}{1+\tn^{1+\varepsilon}})^{-A}.
\end{align}
By \eqref{eq:stirlingmain} we have that
\begin{align}
\frac{\Gamma(\frac{1+\sigma-i\tau+(k+1)/2}{2})\Gamma(\frac{1+\sigma-i\tau+(k-1)/2}{2})}{\Gamma(\frac{-\sigma+i\tau+(k+1)/2}{2})\Gamma(\frac{-\sigma+i\tau+(k-1)/2}{2})} &\ll_\sigma (|\tau|^2 + k^2)^{\sigma+\onehalf}
\end{align}
and hence

\begin{align}
\label{eq:gl2firstbound}
\Psi(x) &\ll_{\sigma,A} \displaystyle\int_{-\infty}^{\infty} (xN)^{-\sigma}(1+\frac{|\tau|}{1+\tn^{1+\varepsilon}})^{-A} (|\tau|^2 + k^2)^{\sigma+\onehalf}d\tau \\ \nonumber
& \ll (1+\tn^{1+\varepsilon})U^{1/2}\left(\frac{U}{xN}\right)^\sigma.
\end{align}

Note that if $xN \ge U(Nk)^\varepsilon$, then taking $\sigma$ large shows that \eqref{eq:gl2firstbound} is consistent with \eqref{eq:gl2boundmain}. So for the rest of the proof we will assume that 
\begin{align}
xN \le U(Nk)^\varepsilon.
\end{align}
Now if $\tn \ll k^\varepsilon$, then we can take $\sigma=0$ to see that \eqref{eq:gl2firstbound} is consistent with \eqref{eq:gl2boundmain}. So we will also assume henceforth that $\tn \gg k^\varepsilon$. For convenience we will also set $\sigma = -\onehalf$. From Lemma \ref{lemma:mellin}, we know that $\tilde{\psi}(-\sigma +i\tau)$ is very small outside of the interval\\$\tn^{1-\varepsilon} \ll |\tau| \ll \tn^{1+\varepsilon}$, so we will restrict integration to this interval. We will now replace $\tilde{\psi}(-\sigma +i\tau)$ in $\Psi(x)$ with the asymptotic formula in \eqref{eq:mellinasymp},
\begin{align}
N^\onehalf W(-\frac{\tau}{\theta})|\tau|^{-\onehalf} e^{i\tau\log |\tau/(e\theta)|} + O(|\tau|^{-3/2}),
\end{align}
where $W$ is a function satisfying \eqref{eq:wconditions}. Define
\begin{align}
\Phi(x) = -\frac{(xN\pi^2)^{1/2}}{2\pi^2}\displaystyle\int_{-\infty}^\infty (x\pi^2)^{i\tau}&W(-\frac{\tau}{\theta})|\tau|^{-\onehalf} e^{i\tau\log |\tau/(e\theta)|}\\\nonumber 
&\frac{\Gamma(\frac{1+\sigma-i\tau+(k+1)/2}{2})\Gamma(\frac{1+\sigma-i\tau+(k-1)/2}{2})}{\Gamma(\frac{-\sigma+i\tau+(k+1)/2}{2})\Gamma(\frac{-\sigma+i\tau+(k-1)/2}{2})} d\tau.
\end{align}
The error term satisfies
\begin{align}
|\Psi(x) - \Phi(x)| \ll \frac{\sqrt{xN}}{(\tn+k)^{100}} + \sqrt{xN}\int_{\tn^{1-\varepsilon} \ll |\tau| \ll \tn^{1+\varepsilon}}|\tau|^{-3/2}d\tau \ll \frac{\sqrt{xN}}{\tn^{1/2-\varepsilon}},
\end{align}
which is satisfactory for \eqref{eq:gl2boundmain}.\\

Using the asymptotic expansion in \eqref{eq:psifirstasymp}, we can write $\Phi(x)$ as a linear combination of expressions of the form $\sqrt{xN}J$ plus an error term, where

\begin{align}
J = \int_{-\infty}^\infty g(\tau)e^{ih(\tau)}d\tau,
\end{align}
where 
\begin{align}
\label{eq:firstphase}
h(\tau) & = \tau\log\left(\frac{(2\pi)^2ex|\tau|}{|\theta|(\tau^2+(k/2)^2)}\right)-k\arctan(\frac{2\tau}{k})
\end{align}
and $g(\tau)$ is a smooth function with support on the interval $|\tau|\asymp\tn$ and satisfying
\begin{align}
\label{eq:gprop}
\frac{d^j}{d\tau^j}g(\tau) \ll |\tau|^{-\onehalf-j}.
\end{align}
Note that the term $\exp(i\arctan(\frac{2\tau}{k}))$ from the asymptotic expansion is considered to be part of the weight function $g$.
Also, note that the error in this asymptotic expansion can be made to be $O(k^{-A})$ for $A$ arbitrarily large, so we only need to bound $J$.

Now we compute some derivatives. Without loss of generality, we will assume that $\tau > 0$.
\begin{align}
\label{eq:firstderiv}
h'(\tau) &= \log\left(\frac{(2\pi)^2x\tau}{|\theta|(\tau^2+(k/2)^2)}\right)\\
\label{eq:holoderiv}
h''(\tau) &= -\frac{1}{\tau}\left(\frac{4\tau^2 - k^2}{4\tau^2+k^2}\right)\\
\label{eq:thirdderiv}
h'''(\tau) &= \frac{16\tau^4-16\tau^2k^2-k^4}{\tau^2(4\tau^2+k^2)^2}.
\end{align}

Notice that $f''(\tau)$ has a zero at $\tau_{00}= k/2$. Since $\tau \asymp \tn$, we will integrate through this zero if $k^{1-\varepsilon} \leq \tn \leq  k^{1+\varepsilon}$. Now when $\tn$ is outside of this range we can apply Lemma 5.1.3 \cite{H} with $V\asymp\tn^{-1/2}$ and $\lambda \gg \tn^{-1+\varepsilon}$ to get $J \ll \tn^\varepsilon$, which is satisfactory for \eqref{eq:gl2boundmain}.

Suppose that $k^{1-\varepsilon} \leq \tn \leq  k^{1+\varepsilon}$. Then when $\tau$ is in a small interval around $\tau_{00}$, say when $|\tau-\tau_{00}| \ll k^{1-\varepsilon}$, the second derivative is too small to use Lemma 5.1.3 \cite{H}. In this case, we write
\begin{align}
h(\tau) &= h(\tau_{00}) + h'(\tau_{00})(\tau-\tau_{00}) + \frac{h'''(\tau_{00})}{6}(\tau-\tau_{00})^3 + H(\tau),
\end{align}
Now in a small interval around $\tau_{00}$, say when $|\tau-\tau_{00}| \ll k^{3/4-\varepsilon}$, $H$ and all its derivatives are small. More precisely $H \ll k^{-\varepsilon}$, and its higher derivatives satisfy
\begin{align}
H' \ll k^{-3/4-\varepsilon}, \quad H'' \ll k^{-3/2-\varepsilon}, \quad H''' \ll k^{-9/4-\varepsilon},
\end{align}
and for $j\ge4$ 
\begin{align}
H^{(j)} = h^{(j)} \ll k^{-j+1}.
\end{align}

Let $w_0$ be a fixed smooth, compactly-supported function, satisfying $w_0(x) = 1$ for $|x|<1$. Then write $J = I_0 + I_1$, where
\begin{align}
I_0 = e^{ih(\tau_{00})}\int_{-\infty}^{\infty}G(\tau)e^{i(h'(\tau_{00})(\tau-\tau_{00})+\frac{h'''(\tau_{00})}{6}(\tau-\tau_{00})^3)}d\tau,
\end{align}
and
\begin{align}
G(\tau) = g(\tau)w_0(\frac{\tau-\tau_{00}}{k^{3/4-\varepsilon}})e^{iH(\tau)}.
\end{align}
Note that $G^{(j)} \ll k^{-1/2}\left(\frac{1}{k^{3/4-\varepsilon}}\right)^j$, and the Fourier transform of $G(\tau)$ satisfies the bound $\hat{G}(y) \ll k^{1/4}(yk^{3/4})^{-A}$ for arbitrary $A>0$. Hence we can write 
\begin{align}
G(\tau) = \int_{-\infty}^{\infty} \hat{G}(y)e(y\tau)dy = \int_{|y|\le k^{-3/4+\varepsilon}}\hat{G}(y)e(y\tau)dy +  O(k^{-A}).
\end{align} 
Substituting this into the integral and recognizing the Airy function that appears, we find that
\begin{align}
I_0 = \frac{2\pi e^{ih(\tau_{00})}}{(3h'''(\tau_{00}))^{1/3}} \int_{|y|\le k^{-3/4+\varepsilon}}\hat{G}(y)e(y\tau_{00})\text{Ai}(2^{1/3}\frac{h'(\tau_{00})+2\pi y}{h'''(\tau_{00})^{1/3}})  dy +  O(k^{-A}).
\end{align}
\subsection{Aside on the Airy Function}
We will briefly discuss some properties of the Airy function. The Airy function is defined by
\begin{align}
\label{eq:airy}
\text{Ai}(x) = \frac{1}{2\pi} \int_{-\infty}^{\infty} \exp(i\frac{t^3}{3} + ixt)dt.
\end{align}
For large arguments, we have the following asymptotics. As $x$ approaches $+\infty$,
\begin{align}
\label{eq:airy1}
\text{Ai}(x)\sim \frac{e^{-\frac{2}{3} x^{3/2}}}{2\sqrt{\pi}x^{1/4}}.
\end{align}
As $x$ approaches $-\infty$,
\begin{align}
\label{eq:airy2}
\text{Ai}(x)\sim \frac{\sin(\frac{2}{3}x^{3/2} +\frac{\pi}{4})}{\sqrt{\pi}x^{1/4}}.
\end{align}
Note that this function is bounded for all $x$.

\subsection{Back to the Proof}
First consider the case when $h'(\tau_{00}) \gg k^{-2/3+\varepsilon}$. Note that by a Taylor expansion,
\begin{align}
h'(\tau_{00}) = \log\left(\frac{(2\pi)^2xN}{\tn k}\right) = \frac{\Delta}{\tn k}(1+o(1)). 
\end{align}
Now since $h'''(\tau_{00})^{-1/3} \asymp k^{2/3}$ we have by \eqref{eq:airy1} and \eqref{eq:airy2} that
\begin{align}
\label{eq:lastbound}
I_0 &\ll h'''(\tau_{00})^{-1/3+1/12}h'(\tau_{00})^{-1/4} \int_{|y|\le k^{-3/4+\varepsilon}}|\hat{G}(y)|dy \\\nonumber
&\ll  k^{1/2}\frac{(\tn k)^{1/4}}{\Delta^{1/4}}k^{-1/2} \ll \frac{k^{1/2+\varepsilon}}{\Delta^{1/4}}.
\end{align}
The condition $h'(\tau_{00}) \gg k^{-2/3+\varepsilon}$ implies that $\Delta \gg \tn k^{1/3+\varepsilon}$. Hence we've obtained the second bound in \eqref{eq:firstcases}.

Now if $h'(\tau_{00}) \ll k^{-2/3+\varepsilon}$, then we bound the Airy function by a constant, so that
\begin{align}
I_0 &\ll h'''(\tau_{00})^{-1/3} \int_{|y|\le k^{-3/4+\varepsilon}}|\hat{G}(y)| dy\\\nonumber
 &\ll k^{2/3-1/2+\varepsilon} = k^{1/6+\varepsilon}.
\end{align}
In this case, $\Delta \ll \tn k^{1/3+\varepsilon}$ and we have our first bound in \eqref{eq:firstcases}.

Now when $ k^{3/4-\varepsilon} \ll |\tau-\tau_{00}| \ll k^{1-\varepsilon}$, we can chop up the interval into dyadic segments of the form $|\tau-\tau_{00}| \in \gamma = [2^jk^{3/4-\varepsilon}, 2^{j+1}k^{3/4-\varepsilon}]$, where $j = 0,1,2,\ldots$. At most we will need $O(\log k)$ such dyadic intervals to cover our entire interval. We can write $I_1 = \sum_\gamma I_\gamma$, where the sum is over each dyadic interval $\gamma$.

Now by the mean-value theorem we have 
\begin{align}
|h'(\tau)- h'(\tau_{00})| = |h''(\xi)(\tau-\tau_{00})|, 
\end{align}
for some $\xi$.
Hence inside one of these dyadic intervals, which we shall denote by $\gamma = [L, 2L]$, we have
\begin{align}
|h'(\tau)- \frac{\Delta}{\tn k}| \asymp L^2 k^{-2}.
\end{align}
Now first assume that $h'(\tau) \gg L^2k^{-2-\varepsilon}$. Then applying Lemma 5.1.2 \cite{H} with $V \asymp \tn^{-1/2}$, we have that $I_\gamma \ll k^{3/2+\varepsilon}L^{-2}$. Now since $L$ is $\gg k^{3/4-\varepsilon}$, we can conclude that at worst, $I_\gamma \ll k^\varepsilon$, which is satisfactory.

On the other hand, if $h'(\tau) \ll L^2k^{-2-\varepsilon}$, then we must have that $\frac{\Delta}{\tn k} \asymp L^2 k^{-2}$. Applying Lemma 5.1.3 \cite{H} with $\lambda \gg Lk^{-2}$ and $V \asymp \tn^{-1/2}$ gives $I_\gamma \ll L^{-1/2}k^{1/2}$. Solving for $L$, we can write this bound as $I_\gamma \ll \frac{k^{1/2+\varepsilon}}{\Delta^{1/4}}$. But this is the same bound as in \eqref{eq:lastbound}. Note that our bound here holds when $\Delta \gg \tn k^{1/2+\varepsilon}$, while the bound in \eqref{eq:lastbound} holds when $\Delta \gg \tn k^{1/3+\varepsilon}$, but this is okay. This completes the proof of the lemma. \qed

\section{Proof of Theorem 1.2}
\noindent The proof of this theorem will be similar to the previous proof, but significantly less complicated. Again, we will actually be proving a smoothed version of our theorem.
\begin{mytheo}
\label{theorem:secondsmoothed}
Let $u(z)$ be a Hecke Maass form as in Theorem \ref{theorem:second} and $w$ be a smooth weight function as in Theorem \ref{theorem:firstsmoothed}. Then for any real number $\alpha$ and any $\varepsilon>0$,
\begin{align}
\displaystyle\sum_{n\ge 1} \lambda_u(n) e(\alpha n)w(n) \ll_\varepsilon N^{1/2 + \varepsilon}\mathfrak{q}_u^{1/4+\varepsilon},
\end{align}
where $\mathfrak{q}_u\asymp T^2$ is the analytic conductor for $u$.
\end{mytheo}
\subsection{$GL(2)$ Voronoi formula}
Let $u(z)$ be a Maass form with expansion \eqref{eq:massexpansion} and having Laplace eigenvalue $\frac{1}{4}+T^2$. Without loss of generality we will assume that $u$ is either even or odd.

Let $\psi(x)$ a smooth function and compact support on the positive reals. Then for $\sigma > -1$ and $\eta\in {0,1}$, define
\begin{align}
\Psi_\eta(x) = \frac{1}{2\pi i} \displaystyle\int_{(\sigma)}(\pi^2x)^{-s}\frac{\Gamma(\frac{1+s+iT+\eta}{2})\Gamma(\frac{1+s-iT+\eta}{2})}{\Gamma(\frac{-s+iT+\eta}{2})\Gamma(\frac{-s-iT+\eta}{2})}\tilde{\psi}(-s)ds.
\end{align}
Then define
\begin{align}
\Psi_+^e(x) &=\frac{1}{2\pi}(\Psi_0(x)+\Psi_1(x)) & \Psi_+^o(x) &=\frac{1}{2\pi}(\Psi_0(x)-\Psi_1(x))\\
\Psi_-^e(x) &=\frac{1}{2\pi}(\Psi_0(x)-\Psi_1(x)) & \Psi_-^o(x) &=\frac{1}{2\pi}(\Psi_0(x)+\Psi_1(x))
\end{align}
\begin{mytheo} \cite[Equations (1.12), (1.15)]{MS2}
\label{theorem:gl2voronoi}
Let $\psi(x)$ be a smooth function with compact support on the positive reals. Let $d,\overline{d},c \in \mz$ with $c\ne0, (c,d) = 1,$ and $d\overline{d}\equiv1 \pmod{c}$. Then if $u$ is even,
\begin{align}
\label{eq:gl2VoronoiMaass}
\displaystyle\sum_{n\ge 1}\lambda_u(n)e(\frac{nd}{c}) \psi(n) &= c\displaystyle\sum_{n\ge 1}\frac{\lambda_u(n)}{n}e(\frac{n\overline{d}}{c}) \Psi_+^e(\frac{n}{c^2})+ c\displaystyle\sum_{n\ge 1}\frac{\lambda_u(n)}{n}e(\frac{-n\overline{d}}{c}) \Psi_-^e(\frac{n}{c^2}),
\end{align}
and if $u$ is odd, then \eqref{eq:gl2VoronoiMaass} holds with $\Psi_\pm^e$ replaced with $\Psi_\pm^o$.
\end{mytheo}

\subsection{Bounding $\Psi_\eta(x)$ and $S$}
Now we let $Q\ge 1$ be a parameter to be chosen later and choose $a, q,$ and $\theta$ as in the proof of Theorem \ref{theorem:firstsmoothed}. Then we write
\begin{align}
S = \displaystyle\sum_{n\ge 1} \lambda_u(n) e(\frac{an}{q}) \psi(n),
\end{align}
where 
\begin{align}
\label{eq:psi3}
\psi(y) = e^{i\theta y}w(y).
\end{align} 

The following bound, due to Iwaniec \cite{Iw3}, will be useful later.
\begin{align}
\label{eq:iwbound}
\sum_{n\le N} |\lambda_u(n)| \ll_\varepsilon N^{1+\varepsilon}T^\varepsilon.
\end{align}

We now state the analogue of Lemma \ref{lemma:gl2bound}.
\begin{mylemma}
\label{lemma:maassbound}
Let $\psi(x)$ be defined by \eqref{eq:psi3} and define
\begin{align}
U = \max(T^2, \tn^2).
\end{align}
Then 
\begin{align}
\Psi_\eta(x) \ll U^{1/2}|NT|^\varepsilon(1+\frac{xN}{U(NT)^\varepsilon})^{-A}.
\end{align}
\end{mylemma}
\noindent {\it Proof of Theorem \ref{theorem:secondsmoothed}.}
Using Lemma \ref{lemma:maassbound}, Theorem \ref{theorem:gl2voronoi}, and \eqref{eq:iwbound} it is easy to see that
\begin{align}
S \ll (QT + Q^{-1}N)(NTq)^\varepsilon.
\end{align}
Choosing $Q = N^{1/2}T^{-1/2}$ gives the bound in Theorem \ref{theorem:secondsmoothed}. \qed\\

\noindent {\it Proof of Lemma \ref{lemma:maassbound}.} Due to similarities with the previous proof, we will only sketch the proof. By Stirling's approximation, we have that
\begin{align}
\frac{\Gamma(\frac{1+\sigma-it+\eta}{2})}{\Gamma(\frac{-\sigma+it+\eta}{2})} = |t/2|^{\sigma+\onehalf}e^{-it\log|t/2e|}(c_0 +\frac{c_1}{|t|}+\cdots+O(\frac{1}{|t|^A}),
\end{align}
where the $c_j$ are constants depending only on $\eta$ and the sign of $\tau$.
Hence
\begin{align}
\frac{\Gamma(\frac{1+\sigma-i\tau+iT+\eta}{2})\Gamma(\frac{1+\sigma-i\tau-iT+\eta}{2})}{\Gamma(\frac{-\sigma+i\tau+iT+\eta}{2})\Gamma(\frac{-\sigma+i\tau-iT+\eta}{2})} &\ll_\sigma [(1+|\tau-T|)(1+|\tau+T|)]^{\sigma+\onehalf} \\ \nonumber
&\ll (|\tau|^2 + T^2)^{\sigma+\onehalf}
\end{align}
and

\begin{align}
\Psi_\eta(x) &\ll_{\sigma,A} \displaystyle\int_{-\infty}^{\infty} (xN)^{-\sigma}(1+\frac{|\tau|}{1+\tn^{1+\varepsilon}})^{-A} (|\tau|^2 + T^2)^{\sigma+\onehalf}d\tau \\ \nonumber
& \ll (1+\tn^{1+\varepsilon})U^{1/2}\left(\frac{U}{xN}\right)^\sigma.
\end{align}
This bound is satisfactory except possibly when $\tn \gg T^\varepsilon$ and $xN \le U(NT)^\varepsilon$. So we now assume that these conditions hold. Set $\sigma = -\onehalf$.

Using the asymptotic formula in \eqref{eq:mellinasymp}, define
\begin{align}
\Phi_\eta(x) = -\frac{(xN\pi^2)^{1/2}}{2\pi^2}\displaystyle\int_{-\infty}^\infty (x\pi^2)^{i\tau}&W(-\frac{\tau}{\theta})|\tau|^{-\onehalf} e^{i\tau\log |\tau/(e\theta)|}\\\nonumber 
&\frac{\Gamma(\frac{1+\sigma-i\tau+iT+\eta}{2})\Gamma(\frac{1+\sigma-i\tau-iT+\eta}{2})}{\Gamma(\frac{-\sigma+i\tau+iT+\eta}{2})\Gamma(\frac{-\sigma+i\tau-iT+\eta}{2})} d\tau,
\end{align}
where $W$ is a function satisfying \eqref{eq:wconditions}. The error term satisfies $|\Psi_\eta(x) - \Phi_\eta(x)| \ll \frac{\sqrt{xN}}{\tn^{1/2-\varepsilon}}$, which is satisfactory. 

We write $\Phi_\eta(x) = \Phi_1(x) + \Phi_2(x)$, where $\Phi_2$ represents the part of the integral where $|\tau\pm T| \le \sqrt{T}$. In this case $\tn \asymp T$ and a trivial bound gives $\Phi_2 \ll \sqrt{xN}$, which is consistent with our desired bound. Now we write $\Phi_1$ as a linear combination of expressions of the form $\sqrt{xN}J$, where

\begin{align}
J = \int_{|\tau\pm T| > \sqrt{T}} g(\tau)e^{ih(\tau)}d\tau,
\end{align}
where 
\begin{align}
h(\tau) & = \tau\log\left(\frac{\pi^2x|\tau|}{e|\theta|}\right) - (\tau+T)\log(|\tau+T|/2e) - (\tau-T)\log(|\tau-t|/2e)
\end{align}
and $g(\tau)$ is a smooth function with support on the interval $|\tau|\asymp\tn$ and satisfying \eqref{eq:gprop}. Note that the error in our expansion can be made to be $O(T^{-A})$ for $A$ arbitrarily large, so we only need to bound $J$.

Without loss of generality, we assume that $\tau > 0$ and compute the derivatives:
\begin{align}
h'(\tau) &= \log\left(\frac{(2\pi)^2x\tau}{|\theta (\tau^2-T^2)|}\right)\\
\label{eq:maassderiv}
h''(\tau) &=-\frac{1}{\tau}\left(\frac{\tau^2 + T^2}{\tau^2-T^2}\right).
\end{align}
Applying Lemma 5.1.3 \cite{H} with $V\asymp\tn^{-1/2}$ and $\lambda \gg \tn^{-1+\varepsilon}$ gives $J \ll \tn^\varepsilon$. Hence we have proven
Lemma \ref{lemma:maassbound}. \qed\\

Notice that this proof was much simpler than the holomorphic case because here $h''(\tau)$ is zero-free. (Compare equation \eqref{eq:maassderiv} with \eqref{eq:holoderiv}).

\section{Proof of Theorem 1.3}
\noindent Again, we prove a smoothed version of our theorem.
\begin{mytheo}
\label{theorem:thirdsmoothed}
Let $F$ be as in Theorem \ref{theorem:first} and let $w$ be as in Theorem \ref{theorem:firstsmoothed}. Then for any real number $\alpha$ and any $\varepsilon>0$,
\begin{align}
\label{eq:thirdsmoothed}
\displaystyle\sum_{n\ge 1} A_F(1,n) e(\alpha n)w(n) \ll_{\varepsilon,c_j} N^{3/4 + \varepsilon}\mathfrak{q}_F^{1/4+\varepsilon},
\end{align}
where $\mathfrak{q}_F\asymp k^2$ is the analytic conductor for $F$.
\end{mytheo}
\subsection{GL(3) Voronoi Formula}
It is known that the symmetric-square lift of a holomorphic modular form is associated with $GL(3,\mz)$-automorphic distribution and so we will use the $GL(3)$ Voronoi summation formula proven by Miller and Schmid \cite[Theorem 1.18]{MS2}. We apply their theorem with the following parameters (see \cite[Proposition 5.12]{MS1}):
\begin{align}
\lambda & = (1-k, k-1, 0)\\
\delta & = (1,0,1).
\end{align}

Let $\psi(x)$ be a smooth function with compact support on the positive reals. Then for $\sigma > -1$ and $\eta \in \{0,1\}$, define 
\begin{align}
\label{eq:psieta}
\Psi_\eta(x) = \frac{1}{2\pi i} \displaystyle\int_{(\sigma)}(\pi^3x)^{-s}\frac{\Gamma(\frac{1+s+k-\eta}{2})\Gamma(\frac{1+s+k-1 +\eta}{2})\Gamma(\frac{1+s+1-\eta}{2})}{\Gamma(\frac{-s+k-\eta}{2})\Gamma(\frac{-s+k-1+\eta}{2})\Gamma(\frac{-s+1-\eta}{2})}\tilde{\psi}(-s)ds.
\end{align}
Then define
\begin{align}
\Psi_+(x) &= \frac{1}{2\pi^{3/2}}(\Psi_0(x)-i\Psi_1(x))\\
\label{eq:psipm}
\Psi_-(x) &= \frac{1}{2\pi^{3/2}}(\Psi_0(x)+i\Psi_1(x)).
\end{align}

\begin{mytheo}\cite[Theorem 1.18]{MS2}
\label{thm:gl3voronoi}
Let $\psi(x)$ be a smooth function with compact support on the positive reals. Let $d,\overline{d},c \in \mz$ with $c\ne0, (c,d) = 1,$ and $d\overline{d}\equiv1 \pmod{c}$. Then 
\begin{align}
\displaystyle\sum_{n\ge 1} A_F(1,n)e\left(\frac{n\overline{d}}{c}\right)\psi(n) &= c\displaystyle\sum_{n_1|c}\displaystyle\sum_{n_2\ge 1}  \frac{A_F(n_2,n_1)}{n_1n_2} S(d,n_2;c/n_1)\Psi_+\left(\frac{n_2n_1^2}{c^3}\right) \\ \nonumber
&+ c\displaystyle\sum_{n_1|c}\displaystyle\sum_{n_2\ge 1}  \frac{A_F(n_2,n_1)}{n_1n_2} S(d,-n_2;c/n_1)\Psi_-\left(\frac{n_2n_1^2}{c^3}\right),
\end{align}
where $S(a,b;c)$ is the usual Kloosterman sum.
\end{mytheo}

As before, directly applying Miller and Schmid's formula will give different gamma factors for $\Psi_\eta$. We have rewritten these factors using \eqref{eq:gammarelation} to resemble the gamma factors of the functional equation \eqref{eq:functionalequation}.

We now note down an asymptotic expansion for $\Psi_\eta$:
\begin{align}
\label{eq:thirdasymptotic}
\Psi_\eta(x) = \frac{1}{2\pi i} \displaystyle\int_{(\sigma)}(\pi^3x)^{-s} & \left( |\tau|\frac{\tau^2+k^2}{2^3}\right)^{\sigma+\onehalf} e^{-i\tau\log(|\tau|\frac{\tau^2+k^2}{(2e)^3})} e^{i(2-2k)\arctan(\frac{\tau}{k})} \tilde{\psi}(-s)ds\\\nonumber
&\times\left(c + Q_1(k,\tau) + Q_2(k,\tau) + \cdots + O(|\tau|^{-A})\right).
\end{align}
where $c$ is some constant and each $Q_j(k,\tau) = O(\max (k,|\tau|)^{-j})$ is a rational function. The constant $c$ and the functions depend only on $\eta$ and the sign of $\tau$.

\subsection{Estimating S}

Now let $Q \ge 1$ be a parameter to be chosen later and choose $a,q,$ and $\theta$ as in the proof of Theorem \ref{theorem:firstsmoothed}. Then we can rewrite the left-hand side of \eqref{eq:thirdsmoothed} as
\begin{align}
S = \displaystyle\sum_{n\ge 1} A_F(1,n) e(\frac{an}{q}) \psi(n),
\end{align}
where 
\begin{align}
\label{eq:psi}
\psi(y) = e^{i\theta y}w(y).
\end{align}

\begin{mylemma}
\label{lemma:estimateOfS}
If $\psi$ is as in \eqref{eq:psi} and $\Psi_\pm$ is as in (\ref{eq:psieta}-\ref{eq:psipm}), then
\begin{align}
|S| \ll q^{3/2+\varepsilon}\max_\pm\max_{n_1|q}\displaystyle\sum_{n\ge 1} n^{-1+\varepsilon} |\Psi_{\pm}(\frac{nn_1^2}{q^3})|.
\end{align}
\end{mylemma}
\noindent{\it Proof.} Applying the Voronoi formula and Weil's bound, we have
\begin{align}
|S| \ll q\max_{\pm}\displaystyle\sum_{n_1|q}\displaystyle\sum_{n_2\ge 1} \frac{|A_F(n_2,n_1)|}{n_2n_1}(q/n_1)^{1/2}d(q) |\Psi_{\pm}(\frac{n_2n_1^2}{q^3})|,
\end{align}
where $d(q)$ is the divisor function.\\
By Deligne's bound, $A_F(n_2,n_1) \ll (n_2n_1)^\varepsilon$, and hence we have that
\begin{align}
|S| \ll q^{3/2+\varepsilon}\max_{\pm}\displaystyle\sum_{n_1|q}\displaystyle\sum_{n_2\ge 1} \frac{n_2^{-1+\varepsilon}}{n_1^{3/2-\varepsilon}} |\Psi_{\pm}(\frac{n_2n_1^2}{q^3})|. 
\end{align}
Taking the max over $n_1$ gives the result. \qed \\

\subsection{Bounding $\Psi_\eta(x)$ and $S$}

\begin{mylemma}
\label{lemma:Gl3bound}
Let $\psi(x)$ be defined by \eqref{eq:psi} and define
\begin{align}
U = \max(k^2, \tn k^2, \tn^3)
\end{align}
and
\begin{align}
\Delta = \left |xN-\frac{1}{(2\pi)^3} \tn k^2 \right |.
\end{align}
Then 
\begin{align}
\Psi_\eta(x) \ll \mathcal{M + E},
\end{align}
where
\begin{align}
\label{eq:gl3boundmain}
\mathcal{M} = \max(k, \tn^{3/2})|Nk|^\varepsilon(1+\frac{xN}{U(Nk)^\varepsilon})^{-A}
\end{align}
and
\begin{align}
\label{eq:gl3boundextra}
\mathcal{E} = 
\begin{cases}
\frac{k^2}{\tn^{1/2}} & \text{if }k^{2/3} \le \tn \le k^{1-\varepsilon} \text{ and } \Delta \ll \tn^3\\
\frac{k^3\tn}{\Delta} & \text{if }k^{2/3} \le \tn \le k^{1-\varepsilon} \text{ and } \tn^3 \ll \Delta \ll \tn k^2\\
\tn k\min(1,\frac{k^2}{\Delta}) & \text{if }k^{\varepsilon} \le \tn \le k^{2/3} \text{ and } \Delta \ll \tn k^2\\
0 & \text{otherwise}.
\end{cases}
\end{align}
\end{mylemma}

\noindent{\it Proof of Theorem \ref{theorem:thirdsmoothed}.}
By Lemmas \ref{lemma:estimateOfS} and \ref{lemma:Gl3bound}, we have that $S \ll S_\mathcal{M} + S_{\mathcal{E}}$. It is easy to see that 
\begin{align}
S_\mathcal{M} \ll q^{3/2}(k + \tn^{3/2})(Nkq)^\varepsilon,
\end{align}
and since $q|\theta| \le 2\pi Q^{-1}$, we have that
\begin{align}
S_\mathcal{M} \ll (Q^{3/2}k + N^{3/2}Q^{-3/2})(NkQ)^\varepsilon.
\end{align}

To bound $S_\mathcal{E}$, we examine several cases. Recall that $x = nn_1^2d^3/q^3$.

\textbf{Case 1.} Suppose that $k^{2/3} \le \tn \le k^{1-\varepsilon}$ and $\Delta \ll \tn^3$. Then we have
\begin{align}
\label{eq:boundingS1}
S_\mathcal{E} \ll & q^{3/2}\frac{k^2}{\tn^{1/2}}\left(\frac{\tn^2}{k^2}\right)(Nkq)^\varepsilon \\\nonumber
\ll & q^{3/2}\tn^{3/2}(Nkq)^\varepsilon
\end{align}

Now suppose that $\tn^3 \ll \Delta \ll \tn k^2$. Further suppose that $Y \le \Delta \le 2Y$. Then we can divide $[Y,2Y]$ into $\ll Y\tn^{-3}$ subintervals of length at most $\tn^3$. We then get $O(Y\tn^{-3})$ instances of bounds of the form $q^{3/2}\frac{k^3\tn}{Y}\left(\frac{\tn^2}{k^2}\right)(Nkq)^\varepsilon$. So we obtain the bound

\begin{align}
\label{eq:boundingS2}
S_\mathcal{E} \ll & q^{3/2}\frac{k^3}{\tn^2}\left(\frac{\tn^2}{k^2}\right)(Nkq)^\varepsilon \\\nonumber
\ll & q^{3/2}k(Nkq)^\varepsilon,
\end{align}
which is no worse than \eqref{eq:boundingS1} for this range of $\tn$.

\textbf{Case 2.} Suppose that $k^\varepsilon \le \tn \le k^{2/3}$ and $\Delta \ll k^2$. Then we have
\begin{align}
S_\mathcal{E} \ll & q^{3/2}\tn k \tn^{-1}(Nkq)^\varepsilon \\\nonumber
\ll & q^{3/2}k(Nkq)^\varepsilon.
\end{align} 
For $k^2\ll \Delta \ll \tn k^2$, we precede as before by dividing into subintervals of length at most $k^2$ to get \eqref{eq:boundingS2} again. Putting all this together we get that
\begin{align}
S \ll (Q^{3/2}k + N^{3/2}Q^{-3/2})(NkQ)^\varepsilon.
\end{align}
Choosing $Q = N^{1/2}k^{-1/3}$ gives \eqref{eq:thirdsmoothed}.\qed\\

\noindent{\it Proof of Lemma \ref{lemma:Gl3bound}}.
Due to similarities with the first two proofs, the first part of the proof will only be sketched. By \eqref{eq:stirlingmain} we have that
\begin{align}
\frac{\Gamma(\frac{1+\sigma-i\tau+k-\eta}{2})\Gamma(\frac{1+\sigma-i\tau+k-1 +\eta}{2})\Gamma(\frac{1+\sigma-i\tau+1-\eta}{2})}{\Gamma(\frac{-\sigma+i\tau+k-\eta}{2})\Gamma(\frac{-\sigma+i\tau+k-1+\eta}{2})\Gamma(\frac{-\sigma+i\tau+1-\eta}{2})} &\ll_\sigma [(1+|\tau|)(k+|\tau|)^2]^{\sigma+\onehalf} \\ \nonumber
&\ll (|\tau|^3 + (1+|\tau|)k^2)^{\sigma+\onehalf}.
\end{align}
Hence we have that
\begin{align}
\label{eq:gl3firstbound}
\Psi_\eta(x) &\ll_{\sigma,A} \displaystyle\int_{-\infty}^{\infty} (xN)^{-\sigma}(1+\frac{|\tau|}{1+\tn^{1+\varepsilon}})^{-A} (|\tau|^3 + (1+|\tau|)k^2)^{\sigma+\onehalf} d\tau\\ \nonumber
& \ll (1+\tn^{1+\varepsilon})U^{1/2}\left(\frac{U}{xN}\right)^\sigma.
\end{align}
These bounds are satisfactory except possibly when $\tn \gg k^\varepsilon$ and $xN\le U(NT)^\varepsilon$. So we now assume that these conditions hold. Set $\sigma = -\onehalf$.

Using the asymptotic formula in \eqref{eq:mellinasymp}, define
\begin{align}
\Phi_\eta(x) = -\frac{(xN\pi^3)^{1/2}}{2\pi}\displaystyle\int_{-\infty}^\infty (x\pi^3)^{i\tau}&W(-\frac{\tau}{\theta})|\tau|^{-\onehalf} e^{i\tau\log |\tau/(e\theta)|}\\\nonumber 
&\frac{\Gamma(\frac{1+\sigma-i\tau+k-\eta}{2})\Gamma(\frac{1+\sigma-i\tau+k-1+\eta}{2})\Gamma(\frac{1+\sigma-i\tau+1-\eta}{2})}{\Gamma(\frac{-\sigma+i\tau+k-\eta}{2})\Gamma(\frac{-\sigma+i\tau+k-1+\eta}{2})\Gamma(\frac{-\sigma+i\tau+1-\eta}{2})} d\tau,
\end{align}
where $W$ is a function satisfying \eqref{eq:wconditions}. The error term satisfies $|\Psi_\eta(x) - \Phi_\eta(x)| \ll \frac{\sqrt{xN}}{\tn^{1/2-\varepsilon}}$, which is satisfactory for \eqref{eq:gl3boundmain}.

Applying the asymptotic expansion in \eqref{eq:thirdasymptotic}, we can expand $\Psi(x)$ as the sum of expressions of the form $\sqrt{xN}J$, where
\begin{align}
J = \int_{-\infty}^\infty g(\tau)e^{ih(\tau)}d\tau,
\end{align}
where 
\begin{align}
h(\tau) &= \tau\log\left(\frac{(2\pi)^3e^2xN}{\tn(\tau^2 + k^2)}\right) - 2k\arctan(\frac{\tau}{k})
\end{align}
and $g(\tau)$ is a smooth function with support on the interval $|\tau|\asymp\tn$ and satisfying \eqref{eq:gprop}. Note that the error in our expansion can be made to be $O(\tn^{-A})$ for $A$ arbitrarily large, so we only need to bound $J$. (Recall that we have assumed that $\tn \gg k^\varepsilon$.)

Now we compute the derivatives:
\begin{align}
h'(\tau) &= \log\left(\frac{(2\pi)^3xN}{\tn(\tau^2 + k^2)}\right)\\
h''(\tau) &= -\frac{2\tau}{\tau^2+k^2}.
\end{align}
\textbf{Case 1.} Suppose that $\tn \ge k^{1-\varepsilon}$. Then applying Lemma 5.1.3 \cite{H} with $V \asymp \tn^{-\onehalf}$ and $\lambda \gg \tn^{-1-\varepsilon}$ gives the bound $J \ll \tn^\varepsilon$, which is consistent with \eqref{eq:gl3boundmain}.\\
\textbf{Case 2.} Suppose that $k^{2/3} \le \tn \le k^{1-\varepsilon}$. In this case, $U = \tn k^2$. Notice that
\begin{align}
h'(\tau) &= \log\left(\frac{(2\pi)^3xN}{\tn k^2}\right) - \log(1+\frac{\tau^2}{k^2})\\\nonumber 
& = \log\left(\frac{(2\pi)^3xN}{\tn k^2}\right) - \frac{\tau^2}{k^2}(1+o(1)).
\end{align}
So unless $xN \asymp \tn k^2$, we have that $|h'(\tau)| \gg 1$ and applying Lemma 5.1.2 \cite{H} gives $J \ll \tn^{-\onehalf}$, which is consistent with \eqref{eq:gl3boundmain}.\\
\textbf{Subcase 2.1.} Suppose that $xN \asymp \tn k^2$ and that
\begin{align}
\label{eq:deltabound}
|\log\left(\frac{(2\pi)^3xN}{\tn k^2}\right)| \le 100\frac{\tn^2}{k^2} .
\end{align}
Then a Taylor expansion shows that $\Delta \ll \tn^3$. Applying Lemma 5.1.3 \cite{H} again with $V$ as before and $\lambda \gg \frac{\tn}{k^2}$ gives $J\ll\frac{k}{\tn}$, which is consistent with \eqref{eq:gl3boundextra}.\\
\textbf{Subcase 2.2.} Suppose that $xN \asymp \tn k^2$ and that \eqref{eq:deltabound} does not hold. Then $\tn^3\ll\Delta\ll\tn k^2$ and
by a Taylor expansion, $|h'(t)| \asymp |\log\left(\frac{(2\pi)^3xN}{\tn k^2}\right)| \asymp \frac{\Delta}{\tn k^2}$. So by Lemma 5.1.2 \cite{H}, we get $J \ll \frac{\tn^{1/2}k^2}{\Delta}$, which is consistent with \eqref{eq:gl3boundextra}.\\
\textbf{Case 3.}
Suppose that $k^{\varepsilon} \le \tn \le k^{2/3}$. If $\Delta \le 100k^2$, then $J$ is not oscillatory so we do not expect to do better than the trivial bound $J \ll \tn^{1/2}$ (note that in this case $xN \asymp \tn k^2$). If $\Delta > 100k^2$, then the bound from Subcase 2.2 applies.\qed\\

\section{Nonlinear Exponential Sums}

\begin{mytheo}
\label{theorem:sungeneral}
Let $g(t)$ be a real-valued function such that for some fixed positive constant $A$, we have that
\begin{align}
\label{eq:cond1}
0 < |g'(t)| \le FN^{-1} < \frac{1}{100}
\end{align}
and
\begin{align}
\label{eq:cond2}
FN^{-2} \le |g''(t)| \le AFN^{-2}
\end{align}
on the interval $[N,2N]$. Then for any $\alpha \in \mr$,
\begin{align}
S = \displaystyle\sum_{N < n\le 2N} \lambda(n) e(g(n) + \alpha n) \ll N^{1/2}(F^{1/2}+\log N), 
\end{align}
where the implied constant depends only on $A$ and $f$.
\end{mytheo}

The proof of this theorem will require two propositions.
\begin{myprop} \cite[Corollary 8.11]{IK}
\label{prop:first}
Let $h(t)$ be a real function with $\nu \le |h'(t)| \le 1-\nu$ and $h''(t) \ne 0$ on $[a,b]$. Then
\begin{align}
\displaystyle\sum_{a<n<b} e(h(n)) \ll \nu^{-1},
\end{align} 
where the implied constant is absolute.
\end{myprop}

\begin{myprop} \cite[Corollary 8.13]{IK}
\label{prop:second}
Let $h(t)$ be a real function with $0 < \Lambda \le h''(t) \le \eta\Lambda$ on $[a,b]$ with $\eta \ge 1$. Then
\begin{align}
\displaystyle\sum_{a<n<b} e(h(n)) \ll \eta\Lambda^{1/2}(b-a) + \Lambda^{-1/2},
\end{align} 
where the implied constant is absolute.
\end{myprop}
\noindent Note that the positivity of $h''(t)$ is actually unneccesary since conjugating the sum does not change the bounds.
\begin{mylemma}
\label{lemma:sun}
Let $g(t)$ be as in \eqref{eq:cond1}-\eqref{eq:cond2}. Then for $x$ positive and real,
\begin{align}
\displaystyle\sum_{N<n\le2N} e(g(n) \pm xn) \ll
\begin{cases}
NF^{-1/2} & \text{ if } 0 \le x \le 10FN^{-1}, \\
x^{-1} & \text{ if } 10FN^{-1} \le x \le \onehalf. \\
\end{cases}
\end{align}
\end{mylemma}
\noindent {\it Proof.} When $0 \le x \le 10FN^{-1}$, Proposition \ref{prop:second} gives us the desired bound. When $10FN^{-1} \le x \le \onehalf$, we have that
\begin{align}
\frac{x}{2} \le |g'(t)\pm x| \le 1 - \frac{x}{2}.
\end{align}
So Proposition \ref{prop:first} gives the desired bound. \qed
\\\\
{\it Proof of Theorem \ref{theorem:sungeneral}.} Since $f$ is a cusp form, the function $F(z) = y^{k/2}|f(z)|$ is bounded on the upper half-plane. Hence for Im$z>0$, 
\begin{align}
\label{eq:fbound}
f(z) \ll y^{-k/2}, 
\end{align}
where the implied constant depends on $f$. Now the Fourier coefficients of $f$ are given by
\begin{align}
\lambda_f(n) = \displaystyle\int_{-\onehalf}^\onehalf n^{\frac{1-k}{2}}f(z)e(-nz)dx.
\end{align}
Writing $z = x+iy$ and changing variables, we have that
\begin{align}
S = \displaystyle\int_{-\onehalf}^\onehalf f(z+\alpha) \displaystyle\sum_{N < n \le 2N} n^{\frac{1-k}{2}}e^{2\pi n y}e(g(n)-xn)dx.
\end{align}
Setting $y=N^{-1}$ and applying Lemma \ref{lemma:sun}, our bound for $f$, and partial summation, we have that
\begin{align}
S &\ll N^{1/2}\left(NF^{-1/2} \displaystyle\int_0^{10FN^{-1}} dx +\displaystyle\int_{10FN^{-1}}^\onehalf \frac{dx}{x}\right)\\\nonumber
  &\ll N^{1/2}(F^{1/2} + \log N).
\end{align} \qed\\
Theorem \ref{theorem:sun} follows immediately from Theorem \ref{theorem:sungeneral}.

\end{document}